\documentclass[conference]{ieeeconf}
\renewcommand{\baselinestretch}{0.998}

\renewcommand{\textfloatsep}{5mm} 

\IEEEoverridecommandlockouts

\usepackage{cite}
\usepackage[acronym]{glossaries}
\usepackage{amsmath,amssymb,amsfonts}
\usepackage{algorithmic}
\usepackage[pdftex]{graphicx}
\usepackage{xcolor}
\def\BibTeX{{\rm B\kern-.05em{\sc i\kern-.025em b}\kern-.08em
    T\kern-.1667em\lower.7ex\hbox{E}\kern-.125emX}}

\usepackage{units}
\usepackage{mathtools} 
\usepackage{amsthm}
\usepackage{lettrine}

\usepackage{scrextend} 

\usepackage{tabularx}
\usepackage[english]{babel}
\usepackage{cite}

\usepackage{array}
\usepackage{color}
\usepackage{amssymb}
\usepackage{multicol}

\DeclareUnicodeCharacter{2212}{-}

\usepackage{todonotes}

\usepackage{tikz}
\usetikzlibrary{arrows,calc}
\usepackage{pgfplots}
\usepgfplotslibrary{groupplots}

\usepackage{amsthm}
\usepackage{algorithm}
\usepackage{hyperref}
\usepackage{epstopdf}

\newcounter{thm}
\newtheoremstyle{mystyle}
  {}
  {}
  {}
  {}
  {\bfseries\color{black}}
  {.}
  {\newline}
  {\thmname{#1}\thmnumber{ #2: }\normalfont\color{black}\thmnote{ \textit{(#3)}}}
\theoremstyle{mystyle}
\newtheorem{prob}[thm]{Problem}

\newtheorem{theorem}{Theorem}[section]

\newacronym{acr:cvt}{CVT}{continuously variable transmission}
\newacronym{acr:dp}{DP}{dynamic programming}
\newacronym{acr:ecms}{ECMS}{equivalent consumption minimization strategies}
\newacronym{acr:eltms}{ELTMS}{equivalent lap time minimization strategies}
\newacronym{acr:em}{EM}{electric motor}
\newacronym{acr:es2k}{ES2K}{Energy Storage to Kinetic}
\newacronym{acr:F1}{F1}{Formula 1}
\newacronym{acr:FIA}{FIA}{F\'{e}d\'{e}ration Internationale de l'Automobile}
\newacronym{acr:fgt}{FGT}{fixed-gear transmission}
\newacronym{acr:FD}{FD}{final drive}
\newacronym{acr:ice}{ICE}{internal combustion engine}
\newacronym{acr:k2es}{K2ES}{Kinetic to Energy Storage}
\newacronym{acr:mgu}{MGU}{motor generator unit}
\newacronym{acr:mguh}{MGU-H}{motor generator unit heat}
\newacronym{acr:mguk}{MGU-K}{motor generator unit kinetic}
\newacronym{acr:mpc}{MPC}{model predictive control}
\newacronym[description={Energy Management Strategy}, \glslongpluralkey={Energy Management Strategies},\glsshortpluralkey={EMSs}]{EMS}{EMS}{Energy Management Strategy}%
\newacronym{acr:pmp}{PMP}{Pontryagin's Minimum Principle}
\newacronym{acr:pu}{PU}{power unit}
\newacronym[description={Powertrain Operation}, \glslongpluralkey={Powertrain Operations},\glsshortpluralkey={POs}]{acr:PO}{PO}{Powertrain Operation}%
\newacronym{acr:rmse}{RMSE}{root-mean-square error}
\newacronym{acr:socp}{SOCP}{second-order cone program}
\newacronym{acr:soe}{SoE}{State of Energy}

\newcommand{\pushright}[1]{\ifmeasuring@#1\else\omit\hfill$\displaystyle#1$\fi\ignorespaces}
\newcommand{\pushleft}[1]{\ifmeasuring@#1\else\omit$\displaystyle#1$\hfill\fi\ignorespaces}
\makeatother
\newif\ifmargincomments 
\margincommentstrue

\newif\ifextendedversion 
\extendedversionfalse

\ifmargincomments

\else

\fi
\begin{document}

\title{\bf Human-in-the-loop Energy and Thermal Management for Electric Racing Cars through Optimization-based Control

}

\author{Erik van den Eshof, Jorn van Kampen, Mauro Salazar%
	\thanks{Erik van den Eshof, Jorn van Kampen, and Mauro Salazar are with the Control Systems Technology section, Department of Mechanical Engineering, Eindhoven University of Technology (TU/e), Eindhoven, 5600 MB, The Netherlands.
	E-mails: {\tt\footnotesize r.c.p.v.d.eshof@student.tue.nl},{\tt\footnotesize j.h.e.v.kampen@tue.nl}, {\tt\footnotesize m.r.u.salazar@tue.nl} 
}}

\maketitle
\begin{abstract}
This paper presents an energy and thermal management system for electric race cars, where we tune a lift-off-throttle signal for the driver in real-time to respect energy budgets and thermal constraints. First, we compute the globally optimal state trajectories in a real-time capable solving time, optimizing a 47-kilometer horizon in 2.5 seconds. Next, for safe operation with a human driver, we simplify it to a maximum-power-or-coast operation in full-throttle regions (straights). Thereby, both the positions from which the vehicle should start coasting and the optimal throttle map are subject to tuning. To this end, we define the coasting sections with a threshold on the co-state trajectory of the kinetic energy from the optimal solution. We devise an online implementable bisection algorithm to tune this threshold and adapt it using PI feedback. Finally, we validate the proposed approach for an electric endurance race car and compare three variants with varying implementation challenges: one re-optimizing and updating the kinetic co-state trajectory online, one applying only the bisection algorithm online, and one relying exclusively on feedback control. Our results show that, under typical racing disturbances, our energy management can achieve stint times ranging from less than 0.056\% to 0.22\% slower compared to offline optimization with a priori knowledge of disturbances, paving the way for on-board implementations and testing.
\end{abstract}

\section{Introduction}
\lettrine[findent=2pt]{\textbf{C}}{limate} change and rapid developments in battery technology are driving a shift in the automotive road car industry towards fully electric vehicles. This trend is also noticeable in motorsport. Formula E has been spearheading this change since 2014 as the major electric racing series. While Formula~E has established itself as one of the most popular racing series in the world, electric racing has since then remained fairly exclusive. This is set to change in the coming years as manufacturers investigate possibilities for electric customer racing and the \textit{Fédération Internationale de l'Automobile} (FIA) is laying foundations for an electric sportscar racing series~\cite{Porsche2023,FIA2023}. Sports racing cars are closely related to road cars, making it one of the cheaper and most popular forms of motorsport. As such, electrification of sportscar racing could promote advances in technological development of electric motorsport and road cars alike~\cite{Borsboom2024}. This electrification is not without its challenges. Electric racing cars suffer from the same problem that deters people from buying electric road cars: the restricted range due to low energy density compared to traditional fuel-powered cars. With recent developments in battery technology paving the way for feasible in-race charging~\cite{InMotion,Kampen2024}, the still relatively long charging time compared to traditional refueling means that maximization of driven distance (i.e., laps) within a stint (period of the race between two pit stops) is of paramount importance. The challenge to retain high average speeds in spite of this makes efficient driving significantly more important in electric cars~\cite{FIAFE2021}. Moreover, the importance of energy management is compounded when multiple re-energizing stops are required, common in long-distance sportscar racing~\cite{Kampen2023a,Kampen2024}. This calls for an energy management strategy that minimizes average lap-time given a number of laps to drive and an energy budget. Against this backdrop, this paper proposes a stint-time-optimal control framework that is computationally efficient enough to run in real-time, robust to disturbances and model errors, and safe in operation alongside a driver.
\setlength{\fboxsep}{0pt}
\setlength{\fboxrule}{1pt}
\begin{figure}[t]
    \centering
    \vspace{2mm} 
    \framebox{\includegraphics[width=\linewidth]{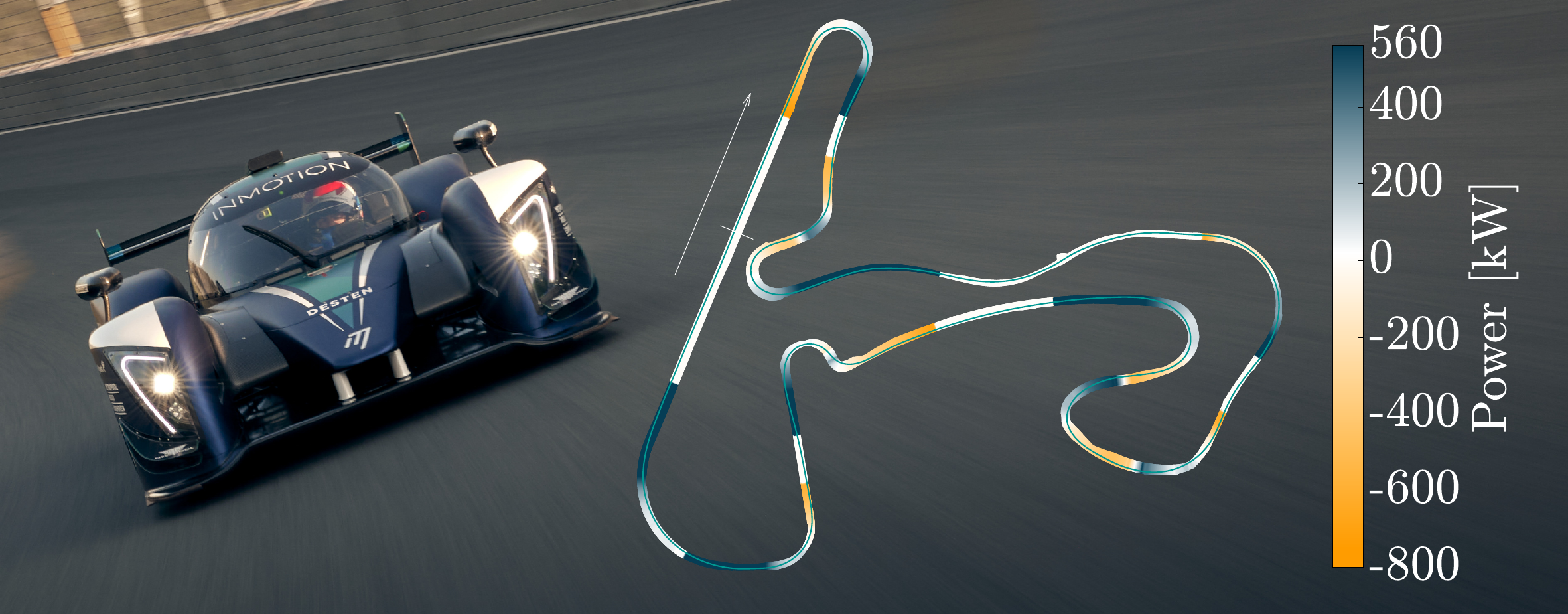}}
    \caption{InMotion's fully electric endurance race car racing around Zandvoort Circuit, with the optimized power trajectory along the circuit to the right.}
    \label{fig:track}
\end{figure}
\paragraph*{\textbf{Related Literature}}
With the introduction of hybrid-electric power-units in Formula 1 in 2014, energy management became an important research stream in motorsport. Rather than direct energy consumption minimization, as in road-cars, the focus in motorsport is on lap-time minimization given, among others, energy consumption constraints~\cite{Borsboom2024}. Researchers have studied racing applications of non-causal techniques such as convex optimization~\cite{Ebbesen2018,Kampen2023a,boyd2004convex}, sensitive to model inaccuracies stemming from convex approximations, and non-linear-programming (NLP)~\cite{Liu2020a}, lacking computational efficiency for real-time application. For real-time optimal control, convex models in combination with Pontryagin's minimum principle (PMP) were studied~\cite{Salazar2017a}, with disturbance rejection through model-predictive control (MPC)~\cite{Salazar2017} and equivalent lap-time minimization strategies (ELTMS)~\cite{Salazar2018}, or a combination of both~\cite{Neumann2023}, but these methods disregard thermal limitations and their application is limited to hybrid powertrains. For hybrid powertrains, the power split between electric and combustion engine is actively controlled and the driver remains in control of the power delivered to the wheels through fixed throttle mappings. This is not the case for fully electric energy management, where optimal control requires an active overwrite of the power delivered to the wheels, which is often restricted by racing regulations and difficult to implement in a safe manner~\cite{F1regs,FEregs}. Borsboom et al. formulated a convex model for lap-time optimal control of electric racing cars, able to solve a horizon of about 14 kilometers in a second \cite{Borsboom2021}, but they disregarded thermal constraints. Van Kampen et al. developed a convex model capable of globally optimal energy and thermal management of an electric race car~\cite{Kampen2023a}, solving a 47 kilometer horizon in 61 seconds, and Herrmann et al. achieved promising results for real-time application using sequential-quadratic-programming (SQP)~\cite{Herrmann2022}. However, their implementations lack proper disturbance rejection and are only suitable for autonomous applications. 
\\
To the best of the author's knowledge, there are no methods specifically focusing on real-time, online implementation of an optimized energy management strategy of electric racing cars. Moreover, thermal constraints are often unaccounted for, and implementation and proof of safe operation alongside drivers is limited.

\paragraph*{\textbf{Statement of Contributions}} This paper presents a robust real-time capable stint-time minimization control framework under energy and thermal constraints. To achieve the computational efficiency necessary for real-time optimal control over an entire stint, we choose a simplified version of the convex model presented and validated by \cite{Kampen2023a}. Subsequently, we adapt the solution to find the optimal lift and coast points and throttle mapping. Moreover, we use feedback control to reduce state errors arising from internal and external disturbances.\\
We quantify the loss of performance of the adaptation method by comparison to the convex solution and we demonstrate the robustness to typical racing disturbances and model errors by comparison to a benchmark optimal solution with disturbances known a priori.



\section{Online Stint Optimization}
\label{ModelSection}
\begin{figure}[t]
    \centering
    \includegraphics[width=0.95\linewidth]{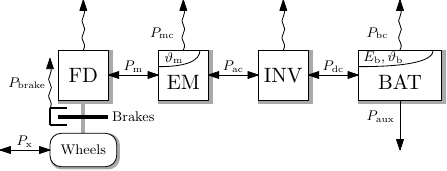}
    \caption{Schematic layout of the electric race car powertrain topology consisting of a battery (BAT), inverter (INV), electric machine (EM), and final drive fixed-gear transmission (FD). The arrows indicate positive power flows of the auxiliary power $P_\mathrm{aux}$, the electrical inverter power $P_\mathrm{dc}$, the electrical EM input power $P_\mathrm{ac}$, the mechanical output power $P_\mathrm{m}$, and the tractive power $P_\mathrm{x}$. Energy losses in the form of dissipated heat are indicated with curvy arrows. The vehicle is rear-wheel driven, with mechanical brakes on each wheel.}
    \label{fig:topology}
\end{figure}

In this section, we illustrate the minimum stint-time control problem capable of real-time, online execution. We leverage the existing convex model presented as low-level model in~\cite{Kampen2023a} and simplify it by condensing the racing line and vehicle-dynamic constraints into a track-dependent maximum kinetic energy constraint~\cite{Salazar2017,Salazar2017a,Salazar2018,Ebbesen2018,Neumann2023}. By this simplification, we only model longitudinal and powertrain dynamics and exclusively focus on propulsion control.
Furthermore, we extend the model with an additional non-convex constraint to model realistic driver behavior and to comply with typical regulatory limits, and provide an algorithm capable of solving this extended problem in real-time as a shrinking horizon model-predictive controller.\\
The electric powertrain topology subject to modeling is shown schematically in Fig.~\ref{fig:topology}. The rear wheels are driven by electric motors (EM) through a fixed gear transmission (FD), powered by inverters (INV) from the battery pack (BAT). We account for a bi-directional energy flow between the wheels and the battery, as the powertrain is capable of charging the battery pack by braking the vehicle (regenerative braking). Additionally, we consider a unidirectional power supply from the battery to unmodeled auxiliary components (Pumps, LV, etc.), and unidirectional energy losses in each of the components in the form of dissipated heat. The EM and battery are actively cooled; therefore, their respective temperatures are also modeled. As we define our model in the space domain, we unconventionally model our energy flows as forces rather than powers. As such, our input variables are the brake force $F_\mathrm{brake}$ and the EM force $F_\mathrm{m}$, both of which represent driver pedal input. As state variables, we use the kinetic energy of the vehicle $E_\mathrm{kin}$, the battery energy $E_\mathrm{b}$, the motor temperature $\vartheta_\mathrm{m}$ and the battery temperature $\vartheta_\mathrm{b}$.

\subsection{Optimization Problem}
\label{OptiSection}
Given a predefined stint length, charging time and maximum kinetic energy profile, we formulate our convex control problem in the space domain with the control variables $\mathbf{u}(s)=\left[F_\mathrm{m}\left(s\right),F_\mathrm{brake}\left(s\right)\right]$ and the state variables $\mathbf{x}(s)=\left[E_\mathrm{kin}\left(s\right),\vartheta_\mathrm{m}\left(s\right),E_\mathrm{b}\left(s\right),\vartheta_\mathrm{b}\left(s\right)\right]$. We minimize the integral of lethargy $\frac{\mathrm{d}t}{\mathrm{d}s}(s)$, the inverse of velocity, in the space domain, which is equivalent to minimizing the total remaining stint time:

\begin{prob}[Convex Stint-Time Optimization]\label{Problem1}
The minimum-stint-time control strategy is the solution of
\newlength{\mylength}
\settowidth{\mylength}{\quad \text{Vehicle Dynamics:}} 

\par\nobreak\vspace{-5pt}
\begingroup
\allowdisplaybreaks
\begin{small}

\begin{equation*}
\begin{aligned}
\min_{u(s),x(s)} \quad & \int_{s_0}^{S_\mathrm{stint}}\frac{\mathrm{d}t}{\mathrm{d}s}(s)\, \mathrm{d}s\\
\makebox[\mylength][l]{\text{s.t.}} & \\
\makebox[\mylength][l]{\quad \text{Electric Machine:}} & \text{\normalfont \cite[eq. (18\;\textendash\;19, 23, 26, 29\;\textendash\;31)]{Kampen2023a}}\\
\makebox[\mylength][l]{\quad \text{Inverter:}} & \text{\normalfont \cite[eq. (34)]{Kampen2023a}}\\
\makebox[\mylength][l]{\quad \text{Battery:}} & \text{\normalfont \cite[eq. (36, 40\;\textendash\;46, 49, 51\;\textendash\;54)]{Kampen2023a}}\\
\makebox[\mylength][l]{\quad} & \text{\normalfont \cite[eq. (56\;\textendash\;58, 60)]{Kampen2023a}}\\
\makebox[\mylength][l]{\quad \text{Vehicle Dynamics:}} & \text{\normalfont \cite[eq. (2\;\textendash\;7, 9\;\textendash\;10, 15)]{Kampen2023a}}\\ 
\makebox[\mylength][l]{\quad}& \;E_\mathrm{kin}(s)\;\,\leq E_\mathrm{kin,max}(s)\\
\makebox[\mylength][l]{\quad \text{Initial Conditions:}} & \;E_\mathrm{kin}(s_0) = \widetilde{E}_\mathrm{kin}(s_0),\;E_\mathrm{b}(s_0) = \widetilde{E}_\mathrm{b}(s_0),\\
\makebox[\mylength][l]{\quad} & \;\vartheta_\mathrm{b}(s_0) \;\;\;= \widetilde{\vartheta}_\mathrm{b}(s_0),\;\;\;\;\!\vartheta_\mathrm{m}(s_0) = \widetilde{\vartheta}_\mathrm{m}(s_0),\\
\end{aligned}
\end{equation*}

\end{small}
\endgroup
\end{prob}
\noindent where $s_0$ is the car's real-time distance into the stint, $S_\mathrm{stint}$ is the total stint distance, $E_\mathrm{kin,max}(s)$ is the predetermined maximum kinetic energy bound and the tilde ($\widetilde{\cdot}$) denotes a measured state.
As the objective and feasible domain are both convex, we can compute a globally optimal solution using standard Non-Linear Programming methods (NLP). Moreover, as the problem is exclusively defined by linear equalities and second-order-conic inequalities (or a subset thereof), we can use specialized Second-Order Conic Programming solvers (SOCP) to improve solving speed.
\begin{figure}[t]
    \centering
    \includegraphics[width=\linewidth]{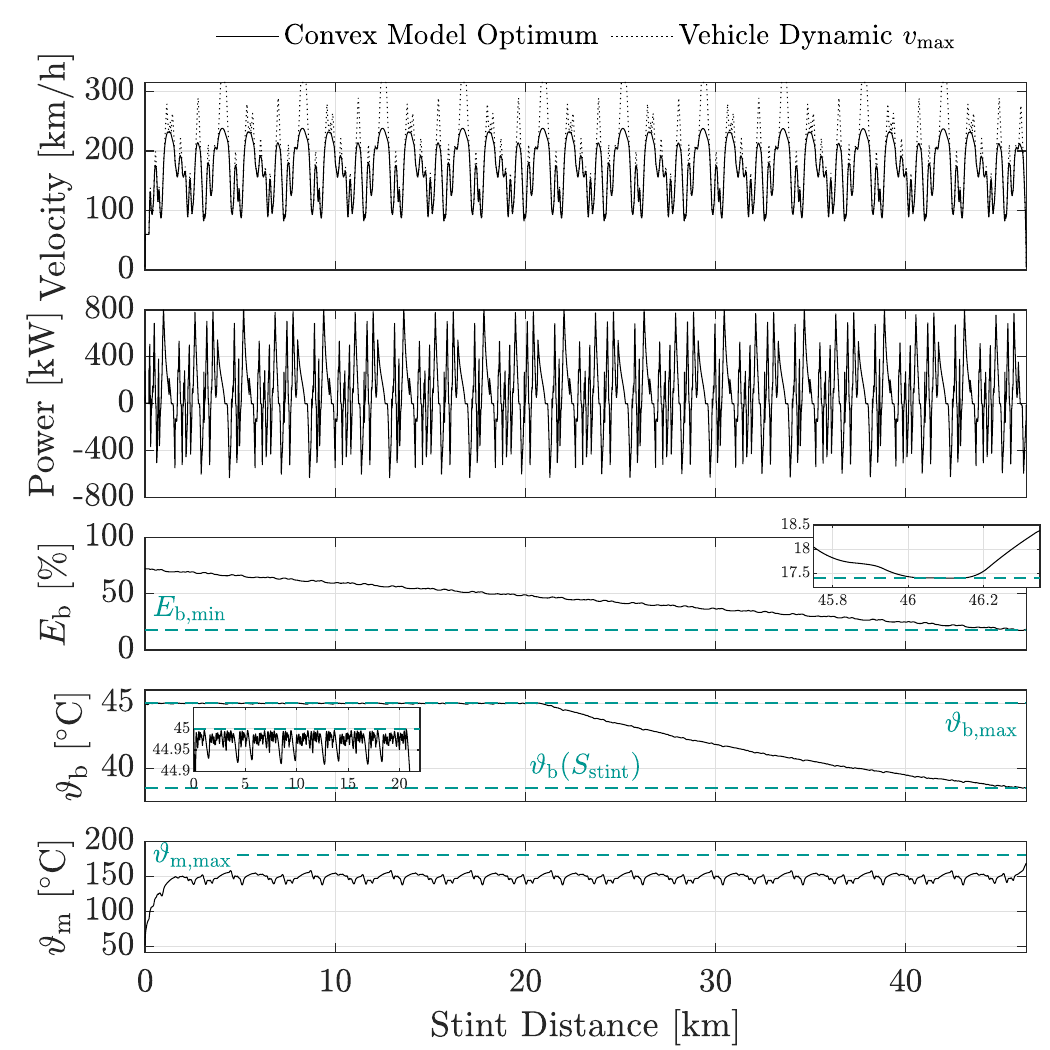}
    \caption{Demonstration of the solution to the convex problem to an 11-lap stint on the Zandvoort circuit. Because of the energy limitation, the power is smoothly reduced on the straights. This smooth trajectory is difficult for a human driver to follow. To maximize efficiency, the battery temperature is initially maximized before cooling down to create margin for fast-charging. The electric machine temperature is not limiting.}
    \label{fig:convexopti}
\end{figure}\\
An example solution to this problem is shown in Fig.~\ref{fig:convexopti}. The solution to this problem is optimal and suitable for autonomous applications. However, it does not guarantee a reference that a human driver can (safely) follow: For safety reasons, racing regulations typically prohibit overwriting of the driver power demand, allowing only fixed \textit{throttle maps} of motor power as a function of throttle percentage and motor rotation speed~\cite{F1regs,FEregs}. Under these regulations the smooth power trajectory shown in Fig.~\ref{fig:convexopti} would require the driver to follow a smoothly varying throttle trajectory when not grip-limited. Therefore, to refrain from enforcing a major distraction and inconvenience on the driver, we include an additional constraint that forces the throttle pedal position $u_\mathrm{th} \in [0,1]$ to be either at $0\%$ or $100\%$ when the car is not limited by grip:
\par\nobreak\vspace{-5pt}
\begingroup
\allowdisplaybreaks
\begin{small}

\begin{equation}
\label{LnCeqn1}
   E_\mathrm{kin}\left(s\right) \neq E_\mathrm{kin,max}\left(s\right) \implies u_\mathrm{th}\left(s\right) = 0 \vee u_\mathrm{th}\left(s\right) = 1.
\end{equation}

\end{small}
\endgroup
As the kinetic energy limit $E_\mathrm{kin,max}(s)$ inherently follows from the driver's ability to consistently drive at- and \textit{sense} the limit of grip, the throttle percentage is free to be varied when grip-limited. This constraint essentially changes the optimization problem to find the optimal points to lift off the throttle pedal to manage energy, which can be communicated to the driver using simple signals like an indicator light or sound~\cite{FormulaE2017}. The throttle pedal position is related to the motor force and vehicle velocity (kinetic energy) through the throttle map(s) $M_i$:
\par\nobreak\vspace{-5pt}
\begingroup
\allowdisplaybreaks
\begin{small}

\begin{equation}
\label{LnCeqn2}
    F_\mathrm{m}(s) = M_i(u_\mathrm{th}(s),E_\mathrm{kin}(s)) \quad \forall i \in \{1,\ldots,n_\mathrm{maps}\},
\end{equation}

\end{small}
\endgroup
\noindent where $n_\mathrm{maps}$ is the number of throttle map options. Using these additional constraints, we reformulate the optimization problem. Inclusion of (\ref{LnCeqn1}--\ref{LnCeqn2}) in Problem~\ref{Problem1} would turn it into a mixed-integer problem, meaning the SOCP solvers that could be used to solve Problem~\ref{Problem1} are no longer applicable. While it is theoretically still possible to find a global optimum using discrete optimization methods like Branch-and-Bound and Dynamic Programming~\cite{Martins2021}, or by turning the problem into a NLP using slack variables~\cite{Duhr2023}, these are not feasible for real-time use due to their long computation times. To still achieve computational efficiency suitable for real-time use, we leverage a two-stage optimization strategy to approximate the solution Problem~\ref{Problem1} s.t. (\ref{LnCeqn1}--\ref{LnCeqn2}). We first solve Problem~\ref{Problem1} using a SOCP solver and extract the kinetic energy co-state $\lambda_\mathrm{kin}(s)\leq0$, which is the Lagrangian multiplier associated with the longitudinal force balance equation ~\cite[eq. (3)]{Kampen2023a} and provided by solvers as one of the outputs. This co-state can be interpreted as a local sensitivity of the objective to changes to the kinetic energy or, in layman's terms, where throughout the stint a decrease in velocity leads to the lowest increase in total stint time. We leverage this information to rewrite (\ref{LnCeqn1}), turning the problem into a 1-dimensional optimization problem for every throttle map where we tune a threshold value $\lambda^*$ to the kinetic co-state that determines where the driver coasts:
\par\nobreak\vspace{-5pt}
\begingroup
\allowdisplaybreaks
\begin{small}
\begin{multline}
\label{LnCeqn3}
\hspace{-3mm}
E_\mathrm{kin}\left(s\right ) \neq E_\mathrm{kin,max} \left(s\right) \implies
\begin{cases}
u_\mathrm{th}\left(s\right) = 0, & \text{if } \lambda_\mathrm{kin}\left(s\right) \geq \lambda_i^*, \\
u_\mathrm{th}\left(s\right) = 1, & \text{if } \lambda_\mathrm{kin}\left(s\right) < \lambda_i^*,
\end{cases}\\
\quad \forall i \in \{1,\ldots,n_\mathrm{maps}\}.
\end{multline}
\end{small}
\endgroup
\begin{figure}[t]
    \centering
    \includegraphics[width=\linewidth]{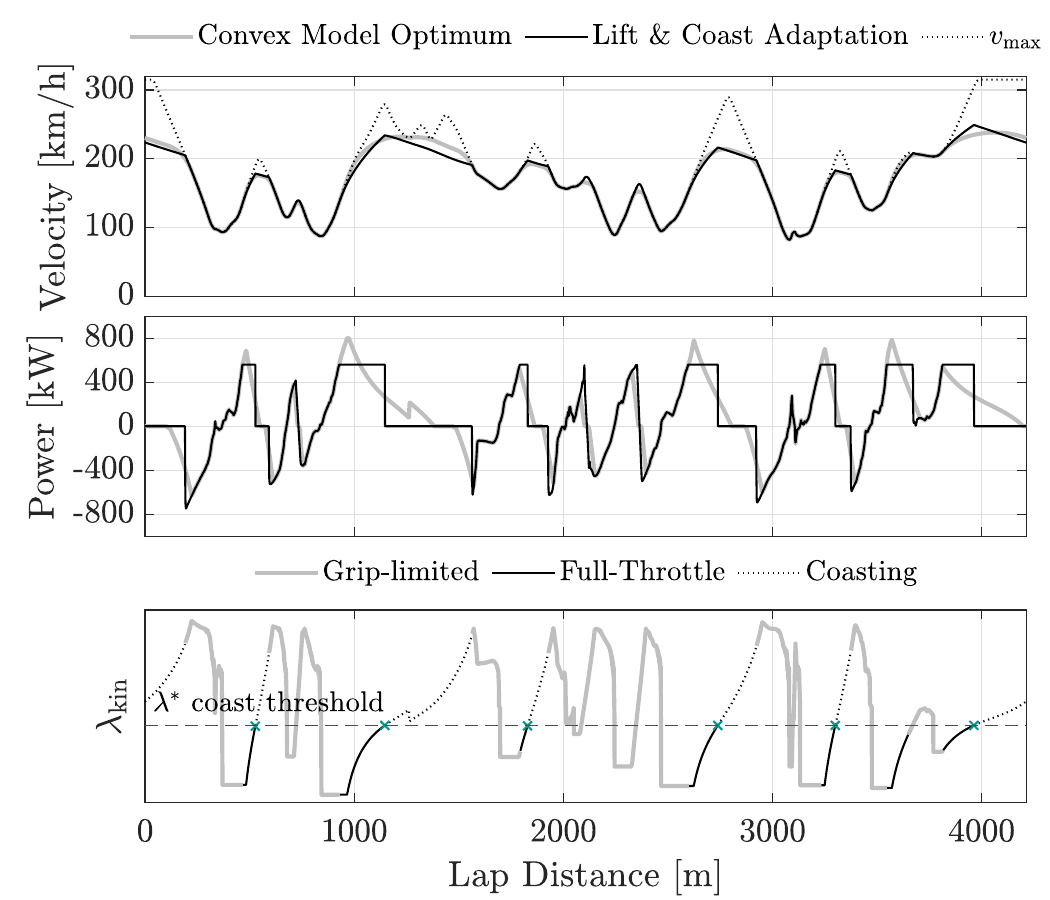}
    \caption{We use the robustified kinetic co-state from the convex problem to enforce the driver constraints as follows: When the driver is full-throttle and the kinetic co-state is above a threshold, a coasting signal is given to the driver. Using the bisection method, the coast threshold is maximized until constraints are violated.}
    \label{fig:lnc}
\end{figure}\\
\noindent This constraint is visually shown in Fig. \ref{fig:lnc}. 
Including these additional constraints, our adapted optimization problem reads as follows:
\begin{prob}[Lift \& Coast Adapted Stint-Time Optimization]\label{Problem2}
The approximated minimum-stint-time control strategy, taking driver constraints into account, is the solution of
\settowidth{\mylength}{\quad \text{\normalfont{Problem~\ref{Problem1} constraints}}} 

\par\nobreak\vspace{-5pt}
\begingroup
\allowdisplaybreaks
\begin{small}

\begin{equation*}
\begin{aligned}
\min_{u(s),\lambda_i^*,x(s),\lambda_\mathrm{kin}(s)} \quad & \int_{s_0}^{S_\mathrm{stint}}\frac{\mathrm{d}t}{\mathrm{d}s}(s)\, \mathrm{d}s\\
\makebox[\mylength][l]{\text{s.t.}} & \\
\makebox[\mylength][l]{\quad \text{Problem~\ref{Problem1} constraints}} & \\
\makebox[\mylength][l]{\quad \text{Drivability:}} & \;(\ref{LnCeqn2}-\ref{LnCeqn3}).
\end{aligned}
\end{equation*}

\end{small}
\endgroup
\end{prob}
\noindent The methodology used to solve this problem is captured in Algorithm~\ref{alg1}. First, Problem~\ref{Problem1} is solved and the trajectory of the kinetic co-state is extracted. Between corner apex (minimum speed) and the co-state minimum of the subsequent straight, we set the co-state to equal this minimum as a robustification measure. Otherwise, when a driver is earlier full-throttle than the algorithm expects, it could occur that a driver is instructed to coast when he is not yet supposed to~\cite{Salazar2018}.
 To optimize the co-state threshold, we use the bisection method~\cite{Vuik2023}. For every threshold value we evaluate, we simulate the complete stint using the same equations governing the convex model, excluding convex relaxations, and apply the second-order Adams-Bashforth integration method~\cite{Vuik2023}. We choose this method because it allows us to analytically invert the model to compute the input required not to exceed a particular state constraint, notably the maximum kinetic energy constraint. To combat the integration error associated with explicit methods, we use a large number of samples, which has no significant influence on the computation time because of the simplicity of the model. We set the motor force to full-throttle or coasting based on the throttle mapping and whether the kinetic co-state is above or below the threshold, and then determine the resulting kinetic energy in the next sampling point. If this kinetic energy exceeds the $E_\mathrm{kin,max}(s)$ bound, we recompute the motor and brake force accordingly so that the kinetic energy is right on the limit instead. As the model evaluated within the bisection loop is not required to be convex, it could optionally be replaced by a higher-fidelity model. However, we chose to use the same model in this study, as in this case the feasible domain of Problem~\ref{Problem1} is a subset of Problem~\ref{Problem2}, meaning that its solution is a lower bound to the achievable cost of Problem~\ref{Problem2}, and thereby a useful benchmark for the algorithm's\label{ResultsLnCAdaptSection}
 performance (see Section~\ref{ResultsLnCAdaptSection}).

 \settowidth{\mylength}{$s_0,E_\mathrm{kin}(s_0),E_\mathrm{b}(s_0),\;$} 
\begin{algorithm}[t]
\caption{Solve Problem~\ref{Problem2}}\label{alg1}
\par\nobreak
\begingroup
\allowdisplaybreaks
\begin{footnotesize}
\begin{algorithmic}
\STATE \makebox[\mylength]{$E_\mathrm{kin,max}(s)
$ \hfill} $\gets \text{predetermined}$\\
\STATE \makebox[\mylength]{$S_\mathrm{stint}$ \hfill} $\gets \text{target pit-stop lap}$\\
\STATE \makebox[\mylength]{$E_\mathrm{b}(S_\mathrm{stint}),\vartheta_\mathrm{b}(S_\mathrm{stint})$ \hfill} $\gets \text{target charging time}$\\
\STATE \makebox[\mylength]{$\begin{array}{l}
s_0,E_\mathrm{kin}(s_0),E_\mathrm{b}(s_0),\;\\
\vartheta_\mathrm{b}(s_0),\vartheta_\mathrm{m}(s_0)
\end{array}$\;\;\;\;\, \hfill} $\gets \text{sensor measurements}$\\
\vspace{1mm}
\STATE $\lambda_\mathrm{kin}(s)$ $\gets \textbf{solve Problem~\ref{Problem1}, robustify}$ \COMMENT{SOCP solver}\\
\vspace{1mm}
\FORALL{$M_i, i \in \{1,\ldots,n_\mathrm{maps}\}$}
\STATE $\lambda_\mathrm{l} \gets \min{(\lambda_\mathrm{kin}(s))}$
\STATE $\lambda_\mathrm{u} \gets \max{(\lambda_\mathrm{kin}(s))}$
\WHILE{$\lambda_\mathrm{u}-\lambda_\mathrm{l} > \text{threshold}$}
\STATE $\lambda_\mathrm{c} \gets (\lambda_\mathrm{u}-\lambda_\mathrm{l})/2$
\STATE $t_\mathrm{stint}, \text{constraints} \gets$ \textbf{evaluate} \textsc{SimulateStint}()
\IF{\text{constraints} \FALSE}
\STATE $\lambda_\mathrm{u} \gets \lambda_\mathrm{c}$
\ELSE
\STATE $\lambda_\mathrm{l} \gets \lambda_\mathrm{c}$
\STATE $\lambda^*_i \gets \lambda_\mathrm{c}$
\STATE $\text{cost}_i \gets t_\mathrm{stint}$
\ENDIF
\ENDWHILE
\ENDFOR
\RETURN $\lambda_\mathrm{kin}(s), \lambda_i^*, \text{cost}_i$
\end{algorithmic}
\vspace{1mm}
{\color{lightgray}\hrule}
\vspace{1mm}
\begin{algorithmic}
\ENSURE{\textsc{SimulateStint}()}
\FOR{$k = s_0$ \TO $k = S_\mathrm{stint}$}
\IF{$\lambda_\mathrm{kin}[k]> \lambda_\mathrm{c}$}
\STATE $F_\mathrm{m}[k] \gets M_i(u_\mathrm{th}=0,E_\mathrm{kin}[k])$ \COMMENT{coast}
\ELSE
\STATE $F_\mathrm{m}[k] \gets M_i(u_\mathrm{th}=1,E_\mathrm{kin}[k])$ \COMMENT{full-throttle}
\ENDIF
\STATE $E_\mathrm{kin}[k+1] \gets \textsc{Model}(F_\mathrm{m}[k])$ \COMMENT{forwards}
\IF{$E_\mathrm{kin}[k+1] > E_\mathrm{kin,max}[k+1]$}
\STATE $E_\mathrm{kin}[k+1]\gets E_\mathrm{kin,max}[k+1]$
\STATE $F_\mathrm{m}[k],F_\mathrm{brake}[k] \gets \textsc{Model}(E_\mathrm{kin}[k+1])$ \COMMENT{backwards}
\ENDIF
\ENDFOR
\RETURN $t_\mathrm{stint}, \text{constraints}$
\end{algorithmic}
\end{footnotesize}
\endgroup
\end{algorithm}

\subsection{MPC Implementation and Feedback}
\begin{figure}[t]
    \centering
    \includegraphics[width=0.95\linewidth]{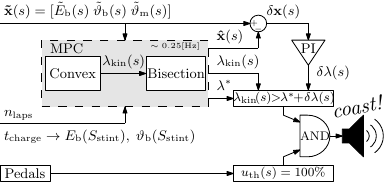}
    \caption{The complete controller framework. Measured signals, target pit stop lap and terminal battery states (dictated by planned charging time $t_\mathrm{charge}$) are the inputs. The MPC controller, limited by solver time, re-optimizes the kinetic co-state trajectory every couple of seconds and provides a threshold. The threshold is adapted using PI controllers based on measured deviations from the optimal state trajectories. If the kinetic co-state is above the threshold and the driver is full-throttle, a signal is given to the driver to start coasting.}
    \label{fig:framework}
\end{figure}

The complete implementation of the framework is shown in Fig.~\ref{fig:framework}.
 As a result of the solving time, the MPC runs at a slower rate than the rest of the framework. Every time the MPC re-initializes, it shrinks its horizon to the remaining distance to be driven and sets its initial states to the measured state values. Based on this, it produces a kinetic co-state trajectory by solving Problem~\ref{Problem1} and a threshold value using the bisection method. We adapt this threshold value by a feedback loop based on the error to the predicted optimal state trajectories using a proportional and integral gain. To prevent integral error wind-up, we reset the integral error every time the MPC updates the trajectories. This feedback loop provides additional resilience to disturbances. Finally, the coasting signal is given to the driver when a number of conditions are met. First and foremost, the kinetic co-state trajectory must exceed the tuned threshold value. Second, as a safety measure and to minimize distraction, the coasting signal can only be given when the driver is full-throttle, i.e., not limited by grip. If necessary, it is possible to add more conditions. For example, a condition could be added that ensures the corner radius is above a certain threshold value, in order to prevent lift-off oversteer (car instability resulting from lifting off-throttle). These conditions should then also be included in the model that is simulated within the bisection loop.

\subsection{Alternative Methods and Discussion}
Due to limitations in the communication between a race car and the garage (regulatory, delays, loss of signal), running the framework on the car is advantageous. However, the complexity of the bisection algorithm and especially the convex solver may be challenging to implement on a race car's electronic control unit (ECU). To investigate the advantages of running these algorithms in real-time on the car, we propose two additional methods that remove the need for running the convex optimizer and the bisection method on the car.
\paragraph*{\textbf{Fully Online}} The full framework is running on the car's ECU, as shown in Fig.~\ref{fig:framework}.
\paragraph*{\textbf{Fixed $\lambda_\mathrm{kin}(s)$}} A fixed kinetic co-state trajectory is assumed. This removes the need to run the convex problem solver on the car's ECU by loading a fixed trajectory on the car, based on a pre-determined optimal strategy. The disadvantage of this method is that it is partially unaware of the vehicle's states, causing a bigger deviation from the optimum depending on the severity of disturbances.
\paragraph*{\textbf{Fixed $\lambda_\mathrm{kin}(s)$, $\lambda^*$}} A fixed kinetic co-state trajectory \textit{and} co-state threshold are assumed. This also removes the need to run the bisection and vehicle model on the car's ECU and eliminates the MPC altogether. This also means no state trajectories are available to base feedback on. Hence, we instead base our feedback on a measured energy usage compared to a target usage, as follows:
\par\nobreak\vspace{-5pt}
\begingroup
\allowdisplaybreaks
\begin{small}
\begin{equation}
    \delta x (s) = \underbrace{\frac{\Delta\bar{E}_\mathrm{kin}+\Delta\bar{E}_\mathrm{b}}{\Delta s}}_\text{measured energy use} - \underbrace{\frac{\widetilde{E}_\mathrm{kin}(s_0) + \widetilde{E}_\mathrm{b}(s_0) - E_\mathrm{bat}(S_\mathrm{stint})}{S_\mathrm{stint}-s_0}}_\text{target energy use},
\end{equation}
\end{small}
\endgroup
\\
\noindent where $\Delta s$ is a window chosen for measurement, e.g., the length of the lap (to minimize fluctuations). Inclusion of the kinetic energy instead of only the battery energy also helps minimize these fluctuations. The additional disadvantage of this method is that it is unable to react to thermal limitations, requiring a more elaborate feedback scheme to affect the co-state threshold based on proximity to temperature limits. However, as thermal limits are far less significant than energy limits, we choose to only control based on the energy limitations.
\begin{figure}[t]
    \centering
    \includegraphics[width=\linewidth]{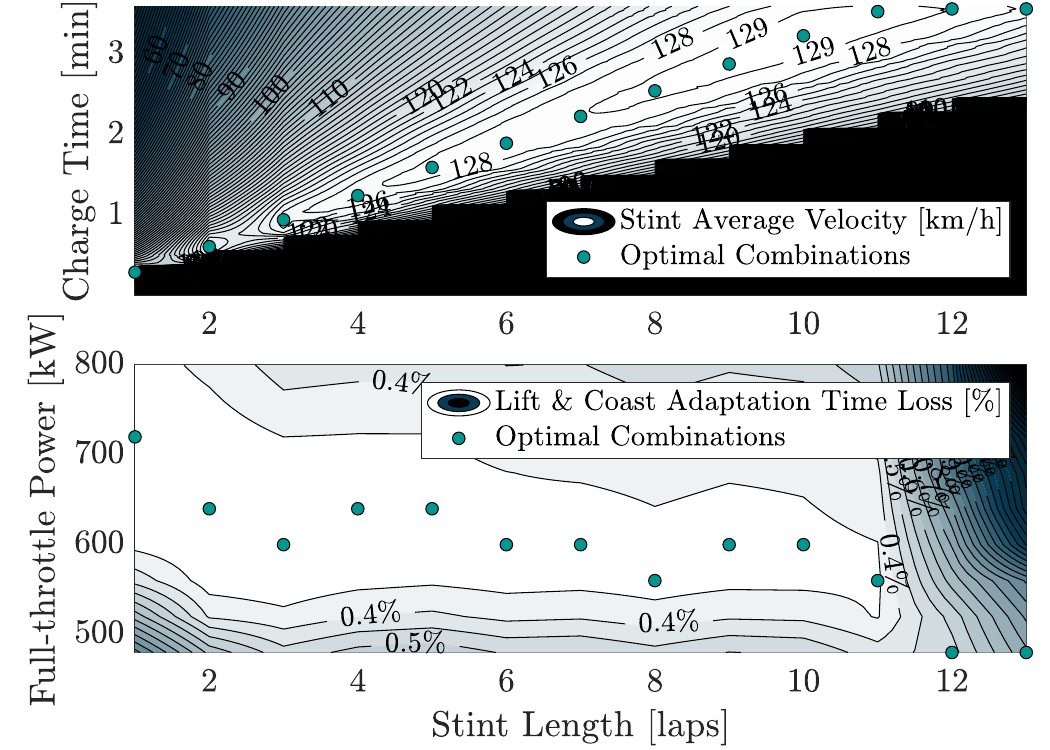}
    \caption{Average stint time (including pit-stop) for a variety of different stint lengths and charging times (pit-stop time). The bottom plot shows, for the optimal combinations of charge time and stint length, the stint-time loss of the lift \& coast adapted solution w.r.t. the convex problem optimum for a variety of different full-throttle powers (throttle maps).}
    \label{fig:lncerror}
\end{figure}
\section{Results}
\label{ResultsSection}
This section presents numerical results on the performance of our MPC framework. First, we validate the performance of the adaptation algorithm by comparing it with the globally optimal solution excluding the full-throttle-or-coast constraint. Next, we look at the results for a full 11-lap stint. Lastly, we subject our framework to various disturbances in real-time simulations and compare its performance against optimization results with \textit{a priori} knowledge of disturbances.\\
We base our results on the electric endurance race car of InMotion~\cite{InMotion}, shown in Fig.~\ref{fig:track}, driving on the Zandvoort circuit. The Zandvoort racing line and maximum kinetic energy profile are precomputed using offline optimization. For the discretization of the convex model running in the MPC, we apply the trapezoidal method with a variable step size ranging from \unit[5]{m} in the corners to \unit[25]{m} on the straights. By this, we can accurately solve stages with quick dynamics while keeping the computation time low. We parse and export Problem~\ref{Problem1} using YALMIP~\cite{Lofberg2004} in MATLAB, and solve the problem using ECOS~\cite{Domahidi2013}, directly updating the sparse input matrices in the MPC to eliminate YALMIP overhead. 
The adaptation algorithm to solve Problem~\ref{Problem2} uses a fixed step size of \unit[1]{m}. We perform the simulations on a standard consumer laptop running Windows with an Intel Core i7-11800H 2.3 GHz processor. Thereby, we achieve consistent solving times of around 2.5 seconds for the SOCP solver and around 0.13 seconds per throttle map for the bisection algorithm, with a horizon of around 47 kilometers.
\begin{figure}[t]
    \centering
    \includegraphics[width=\linewidth]{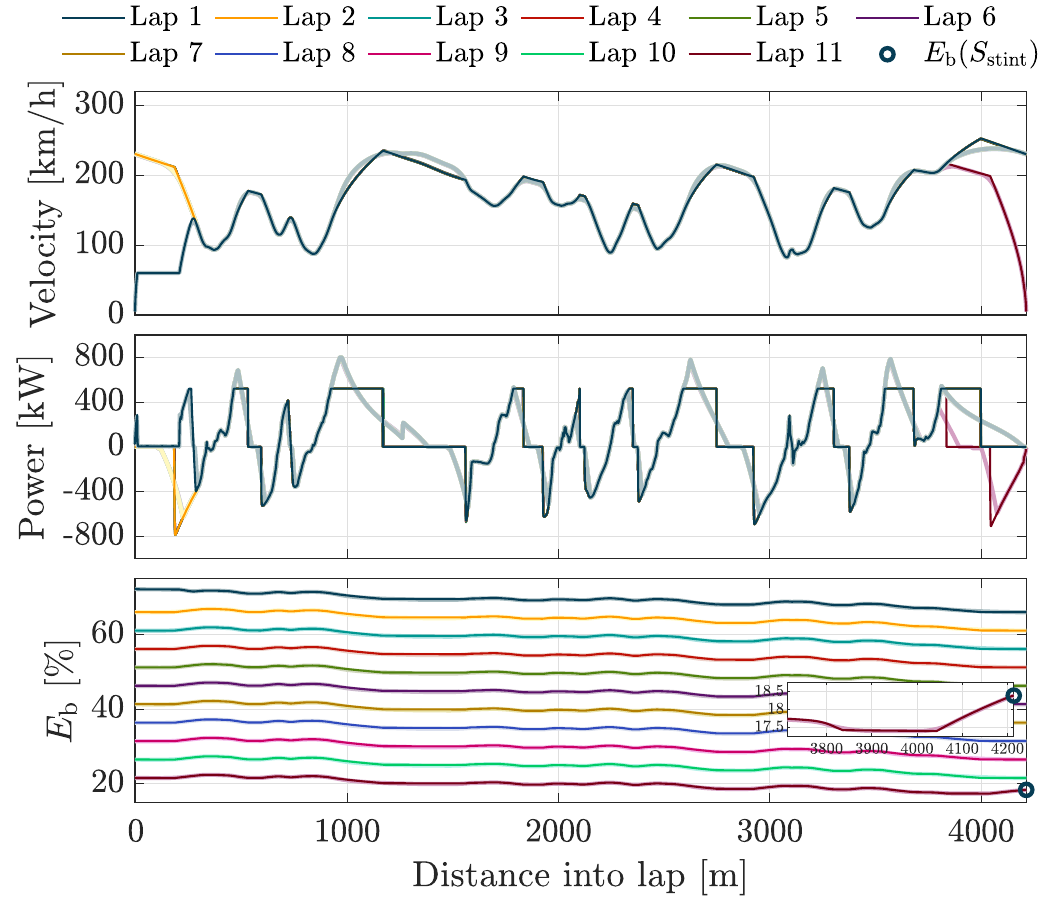}
    \caption{An 11 lap stint with optimized coasting points. The convex model solution is also shown with a lighter color. The coasting points are relatively consistent every lap. The long coasting sections are the result of the significant energy limitation and can be shortened by increasing charging power or artificially limiting the maximum propulsive power.}
    \label{FullStint}
\end{figure}
\subsection{Lift and Coast Adaptation Results}
\label{ResultsLnCAdaptSection}
In this section we discuss the performance and optimality of our lift and coast adaptation algorithm presented in Section~\ref{OptiSection}. First, we determine the optimal combinations of charging time and stint length. Next, we compare the total stint time solution of Problem~\ref{Problem2}, acquired with our adaptation algorithm, to the solution of our convex Problem~\ref{Problem1}. This is done for the optimal combinations of stint length and charging time, and for a variety of throttle maps. For the sake of simplicity, the chosen throttle maps are velocity-independent, meaning that full throttle always corresponds to a fixed maximum power. Different throttle maps may result in a reduction of total stint time, but optimizing the throttle map is beyond the scope of this study.\\
Fig.~\ref{fig:lncerror} shows the average speed during a stint, which includes the objective racing stint time and the additional charging time required to compensate for the energy spent during the stint. The relation between optimal charging time and stint length shows a linear correlation until the capacity of the pack is limiting, for a high number of laps. We also see that, as long as sufficient charging is possible, longer stint lengths result in a higher average velocity because of the lower contribution of the slower pit in- and outlaps. Looking at the time loss compared to the globally optimum solution of  Problem~\ref{Problem1}, we see that this is around \unit[0.3--0.4]{\%}, depending on the throttle map. It is important to note that this is not the loss of our efficient optimization algorithm with respect to the global optimum of Problem~\ref{Problem2}, rather it is an upper bound to the true time loss and thereby an indication of the time lost by enforcing the driver constraints; therefore, we deem this acceptable. 
The small differences between throttle maps indicate that the total stint time is not particularly sensitive to the choice of throttle map, which is advantageous when racing regulations enforce a limited number of throttle maps~\cite{FEregs}.
\\
Fig.~\ref{FullStint} shows a full 11-lap stint, where the terminal battery energy and temperature are based on the planned end-of-stint charging time corresponding to the optimum shown in Fig.~\ref{fig:lncerror}. The initial battery energy and temperature are set to their upper limits, assuming a full, maximum speed charge in the pit-stop before starting the stint. See also Fig.~\ref{fig:convexopti} for the corresponding solution of Problem~\ref{Problem1}. Fig.~\ref{FullStint} reveals that in every lap, the coasting points consistently occur around the same distance into the lap, helping the driver anticipate when they are signaled to start coasting. 
\begin{figure}[t]
    \centering
    \includegraphics[width=\linewidth]{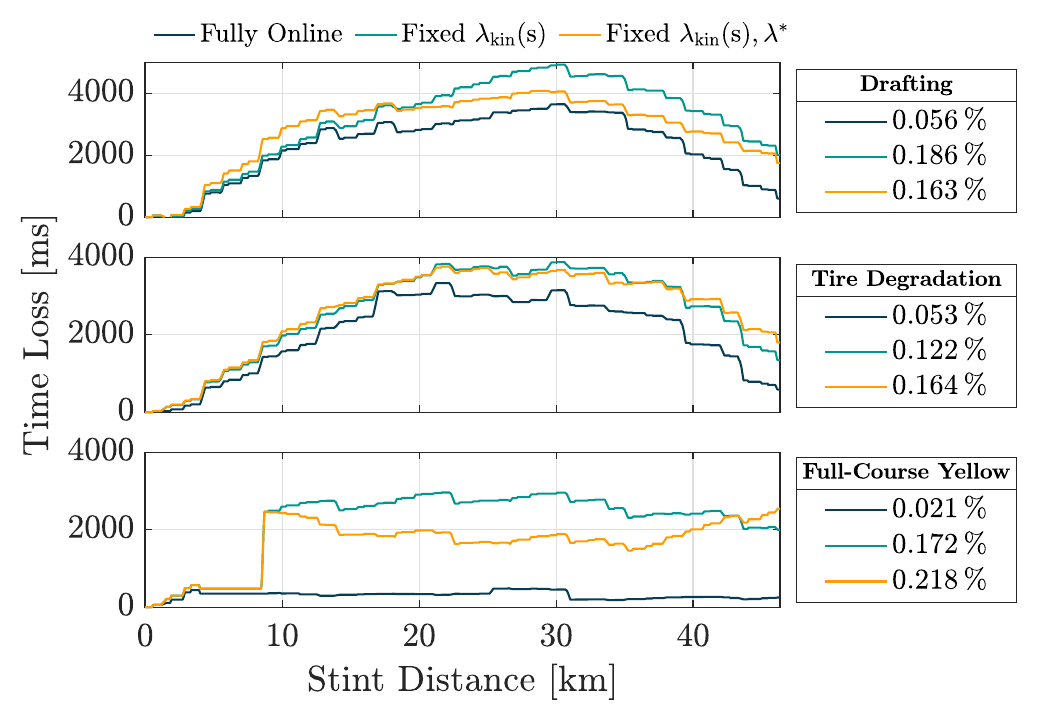}
    \caption{Time loss for real-time simulations of the three algorithms under disturbances to optimization results with a priori knowledge of disturbances. Disturbances are drafting (\unit[90]{\%} down- and dragforce), Tire Degradation (\unit[100]{\%} -- \unit[90]{\%} decreasing grip coefficients over 11 laps), and Full-Course Yellow (2nd lap limited to a maximum speed of \unit[80]{km/h}).}
    \label{fig:distresults}
\end{figure}
\subsection{Results under Disturbances}
In this section, we assess the performance of our methods in 11-lap real-time simulated stints under disturbances, against a non-causal optimization with \textit{a priori} knowledge of the disturbances. We choose three different disturbances that could occur in a real race scenario. First, a drafting disturbance whereby the drag- and downforce are reduced to \unit[90]{\%} of their original value, simulating the effect of following closely behind another car. This affects both the $E_\mathrm{kin,max}(s)$ bound and makes the car more efficient due to reduced drag. Second, a tire degradation scenario where the grip coefficients drop linearly to \unit[90]{\%} of their original values throughout the stint, representing tire wear. Finally a full-course yellow scenario, where the second lap of the stint is restricted to \unit[80]{km/h}, a common safety intervention in racing.\\
The results are shown in Fig.~\ref{fig:distresults}. The fully online approach performs best in all cases, as it is most adaptable to disturbances. The methods with a fixed co-state trajectory or threshold perform comparably; worse than the fully online approach, but not significantly, as typical time deviations in sportscar racing are of much greater magnitude~\cite{Kampen2024,FIAWECTELEMETRY}. The full-course yellow scenario exposes some flaws in the approaches with fixed co-state trajectory: as the full-course yellow ends and driver is accelerating back to racing speed, they are instructed to coast at low speed, as the coasting signal is given at a distance that is comparable to future laps because of the fixed trajectory. Coasting at low speed means significant time loss with little energy saved to make up for it later in the stint. However, a heuristic condition, such as a minimum coasting speed, could prevent this from happening or reduce the significance. For the fixed co-state trajectory threshold, we also see that as a result of integral error wind-up during the full-course yellow period, it starts reducing its time deficit by using more energy, but is then forced to use less energy later in the stint, losing time again. Therefore, it could benefit from an anti-windup measure when coasting is not possible for extended periods. For example, by blocking integral wind-up during full-course yellow periods. As these problems are solvable in a variety of different ways, we leave the determination of the best method to future research and real-world tests.
\section{Conclusion} \label{Conclusion}
In this paper, we presented a framework for online energy and thermal management for electric race cars. First, we devised a method for real-time capable optimization of lift and coast points by leveraging the kinetic energy co-state from a convex optimization solution and tuning a threshold using the bisection method. We showed that this method of enforcing a maximum-power-or-coast constraint only loses around 0.35\% of stint time compared to a globally optimum solution where power is free to vary. Moreover, we showed that the choice of throttle map does not significantly affect the performance of the algorithm, and that in an undisturbed situation the coasting points remain consistent lap-by-lap. Second, we evaluated three approaches with varying levels of challenge to implement on a car. 
We subjected these methods to typical disturbances found in racing in real-time simulations, and compared their time loss with respect to an offline optimization result with a priori knowledge of said disturbances. The method running the complete framework online performed the best, losing around \unit[0.05]{\%} of stint time. The methods that do not actively update the kinetic energy co-state trajectory perform similarly, losing \unit[0.1--0.2]{\%} of stint time. The latter suffers when there is a large deviation in the velocity state, possibly calling to coast at low speed. Moreover, the easiest-to-implement method is unaware of temperature states and is therefore only capable of energy management. However, both of these problems can be improved through heuristic methods and conditions, which should be applied and tested after on-board implementation. Overall, the time loss is insignificant compared to typical deviations in sportscar racing.
Further research and testing on a physical car should refine these methods and determine the most suitable for onboard application. Additionally, the computational efficiency of the lift-and-coast based stint optimization opens the door for a variety of other applications: Higher level strategy optimization (optimizing stint lengths and energy budgets), efficient parameter sensitivity analysis, and usage to warm-start higher-fidelity optimization models.

\section*{Acknowledgment}
\noindent We thank Dr. I. New for proofreading this paper.\\

\bibliographystyle{IEEEtran}
\renewcommand{\baselinestretch}{0.96}

\end{document}